\documentclass[notitlepage,pra,longbibliography,superscriptaddress]{revtex4-1}

\usepackage{amsmath} 
\usepackage{amssymb}
\usepackage{bm}
\usepackage{bbm}
\usepackage{graphicx}
\usepackage{cancel}
\usepackage{xcolor}
\usepackage{tcolorbox}
\usepackage{xfrac}

\usepackage[fleqn]{mathtools}

\DeclarePairedDelimiter\floor{\lfloor}{\rfloor}

\newcommand{\dd}{\mathrm{d}}

\newcommand{\calP}{\mathcal{P}}
\newcommand{\calT}{\mathcal{T}}
\newcommand{\calO}{\mathcal{O}}

\newcommand{\half}{{(\sfrac{1}{2})}}
\newcommand{\expo}{{(1/s)}}

\definecolor{light}{gray}{0.90}
\definecolor{darker}{gray}{0.50}
\definecolor{dark}{gray}{0.30}

\definecolor{garrosgreen}{rgb}{0.1, 0.4, 0.1}
\definecolor{dartmouthgreen}{rgb}{0.05, 0.5, 0.06}
\definecolor{duelferred}{rgb}{0.7, 0.2, 0.1}
\definecolor{cambridgeblue}{rgb}{0.1, 0.3, 1.0}
\definecolor{oxfordblue}{rgb}{0.05, 0.2, 0.7}

\begin{document}

\title{Enhanced and Generalized One--Step Neville Algorithm:\\
Fractional Powers and Access to the Convergence Rate}

\author{Ulrich D.~Jentschura}

\affiliation{Department of Physics and LAMOR, Missouri University of Science and
Technology, Rolla, Missouri 65409, USA}

\author{Ludovico T.~Giorgini}

\affiliation{Nordita, Royal Institute of Technology and Stockholm University,
Stockholm 106 91, Sweden}

\begin{abstract}

The recursive Neville algorithm allows one to 
calculate interpolating functions recursively.
Upon a judicious choice of the abscissas 
used for the interpolation (and extrapolation), 
this algorithm leads to a method for convergence 
acceleration. For example, one can use the 
Neville algorithm in order to successively eliminate 
inverse powers of the upper limit of the 
summation from the partial sums of a given,
slowly convergent input series.
Here, we show that, for a particular 
choice of the abscissas used for the 
extrapolation, one can replace the recursive Neville scheme
by a simple one-step transformation,
while also obtaining access to subleading 
terms for the transformed series after
convergence acceleration.
The matrix-based, unified formulas 
allow one to estimate the rate of convergence 
of the partial sums of the input series
to their limit. In particular, 
Bethe logarithms for hydrogen are calculated 
to 100 decimal digits.
Generalizations of the 
method to series whose remainder terms can be 
expanded in terms of inverse factorial series,
or series with half-integer powers, are also 
discussed.
\end{abstract}

\maketitle

\tableofcontents

%
%
\section{Introduction}
\label{sec1}

Often in physics, one faces problems connected with slowly
convergent, nonalternating series.  Examples include angular momentum
expansions in both nonrelativistic as well as relativistic atomic physics
calculations~\cite{JeMoSo1999,JeMoSoWe1999}, series representations of Bethe
logarithms~(Chap.~4 of Ref.~\cite{JeAd2022book}), or, logarithmic sums over
eigenvalues required for the calculation of the functional determinant in
$O(N)$ theories~\cite{GiEtAl2020,GiEtAl2022,GiEtAl2024i}. The
acceleration of the convergence
these series is notoriously problematic~\cite{We1989,BrRZ1991,CaEtAl2007}. 
Traditional methods like the
Aitken $\Delta^2$ process~\cite{Ai1926}, the Shanks
transformation~\cite{Sh1955}, and Wynn's epsilon algorithm~\cite{Wy1956eps}
(which calculates Pad\'{e} approximants), and other methods, all encounter
numerical instabilities.

For nonalternating series, the challenge of the convergence 
acceleration differs significantly from the resummation
of a divergent input series. While both scenarios aim to
derive meaningful information from the series, the nature of the series and the
techniques applied are fundamentally distinct. The latter seeks to assign a
finite sum to a series that lacks natural convergence. This process often
involves the employment of analytical continuations or summation techniques
that go beyond conventional convergence acceleration methods.
For example, in Ref.~\cite{Je2000prd}, it is shown how, by a suitable
resummation procedure, a dissipative effect, namely, a nonperturbative pair 
creation amplitude in a strong electric background field, 
can be derived from a nonalternating series that sums the 
higher-order dispersive terms of the quantum corrections
to the Maxwell Lagrangian. These quantum electrodynamic corrections
are summarized in a nonalternating divergent perturbative 
expansion of the Heisenberg--Euler Lagrangian
in powers of the electric field strength~\cite{HeEu1936,Sc1951,DiRe1985,JeAd2022book}.
In contrast, summing a slowly convergent series entails enhancing the rate at
which an already convergent series approaches its limit.

For alternating series, it is known that 
nonlinear sequence transformations such as the ones
described in Ref.~\cite{We1989}, notably,
the Weniger--Levin transformation
(see Refs.~\cite{Si1979,Si1980,Si1986,We1989,We1996c,Si2003}), can lead both 
to the acceleration of convergence and also,
to efficient convergence acceleration methods.
For slowly convergent nonalternating series, 
the methods described for their resummation (see Ref.~\cite{Je2000prd})
typically differ drastically from those employed for the 
acceleration of convergence.
In Refs.~\cite{NP1961,JeMoSoWe1999,AkSaJeBeSoMo2003}, an algorithm
was described which first converts 
the nonalternating input series
to an alternating series,
whose convergence is accelerated via nonlinear 
sequence transformations~\cite{We1989,JeMoSoWe1999,AkSaJeBeSoMo2003}.
In fact, the first step of the two-step condensation,
as described in Ref.~\cite{JeMoSoWe1999},
consists of the van Wijngaarden transformation~\cite{NP1961},
which rearranges the nonalternating 
input series into an alternating series, 
but it requires the calculation 
of individual terms of the input series of very high order
in the summation indices.
It is followed by a nonlinear sequence 
transformations~\cite{Si1979,Si1980,Si1986,We1989,We1996c,Si2003},
leading to what is known as the combined
nonlinear-condensation transformation.
Sometimes, however, the option of the 
calculation of terms of very high order 
in the summation indices is not available.
This problem could occur in various contexts,
e.g., because, on lattices, 
we cannot calculate eigenvalues with 
a high principal quantum number~\cite{GiEtAl2020,GiEtAl2022,GiEtAl2024i}.
In the following, we will not discuss the combined
nonlinear-condensation transformation~\cite{JeMoSoWe1999} 
any further, as it contradicts the main 
thrust of the method described in the current paper,
which pertains to convergence acceleration based on 
a limited number of terms of the input series.


In these cases, established methods like Aitken $\Delta^2$
process~\cite{Ai1926} and Wynn's Epsilon Algorithm~\cite{Wy1956eps} are often
used.  One notes, according to Sec.~2.2.7 of Ref.~\cite{CaEtAl2007}, that the
Aitken $\Delta^2$ process and Wynn's Epsilon Algorithm are special cases of the
general E-algorithm~\cite{Br1980E,Ha1979E}.  Of course, the Aitken $\Delta^2$
process and Wynn's Epsilon Algorithm constitutes viable and widely used
convergence acceleration methods, but numerical instability can be encountered
in higher orders of the transformation.

In view of the numerical instabilities encountered
in the convergence acceleration of these series,
it is desirable to use asymptotic information
about the input series to the maximum extent
possible, in order to achieve the maximum
convergence acceleration already in low 
orders of the transformation. One possibility to achieve
this rapid convergence is to use the 
Neville algorithm, which, 
\emph{a priori},
in its most general form, 
constitutes an interpolation algorithm~\cite{Ne1934}.
It can be used for convergence acceleration,
which is tantamount to extrapolation, 
if one uses the algorithm in order to calculate 
the interpolating polynomial outside 
of the domain of interpolation. One first 
interprets the partial sums of the input series
as values of a function at abscissas which 
tend to a limiting value (say, zero) as the number of 
terms of the input series is increased. 
One possibility is to use $x = 1/(i+1)$ for the 
partial sums $s_i$ of the input series ($i=0,1,2,\dots$).
One then calculates the value of the interpolating
function at the limiting value of the abscissa
(in our example case, at $x=0$) which 
is outside of the domain $x \in (1/(n+1), 1)$ 
used for the interpolation (here, $n$ is the 
maximum summation index used for the input series
in a given order $n$ of the transformation).
The Neville algorithm uses information 
from all partial sums $(s_0, \dots, s_n)$ and hence,
about the asymptotic structure of the input series,
in order to construct the interpolating polynomial in $x$
and can therefore lead to highly efficient
convergence acceleration. 
Its usual formulation is based on a recursive 
three-term recursive scheme [see Eq.~\eqref{three_rec}].
Here, we aim to enhance the Neville algorithm
by replacing the three-term recurrence relation with 
unifying analytic formulas which also give
access to subleading asymptotics of the 
behavior of the partial sums $s_n$ of the 
input series for large $n$.

This paper is organized as follows.
In Sec.~\ref{sec2}, we describe our formulation
of the enhanced Neville algorithm.
The idea behind our formulation is discussed 
in Sec.~\ref{sec2A}, while a comparison to the 
Neville algorithm is presented in Sec.~\ref{sec2B}.
We find universal formulas which allow access to the 
subleading terms for high $n$, in Sec.~\ref{sec2C}.
The performance of the algorithm is compared to 
other established methods in Sec.~\ref{sec3}.
Numerical examples are discussed in Sec.~\ref{sec4},
generalizations of the algorithm
are discussed in Sec.~\ref{sec5},
and conclusions are reserved for Sec.~\ref{sec6}.
Appendix~\ref{appa} lists some higher-order 
coefficients for the subleading asymptotic behavior 
of the input series, and Appendix~\ref{appb}
gives numerical values for hydrogen Bethe logarithms obtained using 
our method.

%
%
\section{Enhanced Neville Algorithm}
\label{sec2}

%
%
\subsection{Formulation of the Algorithm}
\label{sec2A}

One starts with the following relations
for the partial sums $s_n$ of the input series,
and the remainder terms $r_n$,
\begin{equation}
\label{remainder_term}
s_n = \sum_{k=0}^n a_k = s_\infty + r_n \,, \qquad
r_n = - \sum_{k=n+1}^\infty a_k \,.
\end{equation}
(Here, we start the terms from index $k=0$, 
by convention~\cite{We1989}.)

Let us assume that the terms of the infinite series
possess an asymptotic expansion in 
inverse integer powers of the summation index,
\begin{equation}
\label{asymptotic}
a_k = \frac{A}{k^2} + \frac{B}{k^3} + \calO(k^{-4}) \,.
\end{equation}
The idea is that, if $n$ is sufficiently large, we 
can use the asymptotic approximation for the remainder term,
\begin{equation}
r_n = - \sum_{k=n+1}^\infty a_k 
= - \sum_{k=n+1}^\infty \left[ \frac{A}{k^2} + 
\frac{B}{k^3}  + \frac{C}{k^4} + \calO(k^{-5}) \right] \,.
\end{equation}
One can use the relation
\begin{equation}
\sum_{k = n+1}^\infty \frac{1}{k^a} = 
\frac{ (-1)^a }{(a-1)!} \, 
\left. \frac{\partial^{a-1} \psi(z)}{\partial z^{a-1}}  
\right|_{z = n + 1} \,,
\end{equation}
where $\psi(z)$ is the logarithmic derivative of the 
Gamma function, to expand the remainder 
term $r_n$ in inverse power of $n$.
The first terms in this expansion read as follows,
\begin{equation}
\label{remainder}
r_n = \frac{A}{n} + \frac{B-A}{2 n^2} 
+ \frac{A - 3B + 2C}{6 n^3} + \calO(n^{-4}) \,.
\end{equation}
One then tries to eliminate the asymptotic terms
of the remainder function, which are powers of 
$1/n$. To this end, one interprets the 
first $(n+1)$ partial sums $s_i$ with $i \in (0, \ldots, n)$ 
as approximations to the values of a function  $f(x)$
at the abscissa values $x = x_i = 1/i$,
\begin{equation}
s_i \approx f(x = x_i = 1/(i+1)) \,,
\qquad i \in (0, \ldots, n) \,,
\end{equation}
so that $\lim_{\epsilon \to 0^+} f(\epsilon) = s_\infty$.
If we assume that $f(x)$ can be expanded as a power series,
then we can write
\begin{equation}
\label{ansatz_poly}
f(x_i) = \sum_{j=0}^{\infty} c_j \, (x_i)^j 
\approx \sum_{j=0}^{n} c_j \, (x_i)^j
= c_0 + \frac{c_1}{(i+1)} + \frac{c_2}{(i+1)^2} + 
\ldots + \frac{c_n}{(i+1)^n} \,,
\end{equation}
with $n+1$ unknown coefficients $c_j$
($j=0,\dots,n$). 
Given $n+1$ partial sums $s_i$ of the input series
$\{ s_i \}_{n=0}^\infty$, the idea of the algorithm is 
to solve the system of equations
\begin{equation}
s_i = \sum_{j=0}^n c_j(n) \,
\left( \frac{1}{i+1} \right)^j \,,
\qquad i, j \in (0, \ldots, n) \,,
\end{equation}
for the coefficients $c_j(n)$. One then assumes that the 
value of the infinite series is recovered as
$n$ is increased,
\begin{equation}
\label{conv_poly}
s = s_\infty = \lim_{n \to \infty} c_0(n) \,,
\end{equation}
where $c_0(n)$ is the zeroth-order coefficient
of the interpolating polynomial 
$\sum_{j=0}^n c_j(n) \, x^j$, which, at the 
abscissas $x = x_i = 1/(i+1)$, coincides with the 
partial sums $s_i$ of the input series.
In view of the fact that $c_0(n)$ is identical 
to the value of the interpolating polynomial 
at argument $x=0$ while all abscissas $x_i$ are 
positive, the latter step constitutes an 
extrapolation which can lead to convergence 
acceleration~\cite{CaEtAl2007}.

To conclude this section, we observe that,
in the sense of Eq.~\eqref{remainder_term},
the {\em ansatz} given in Eq.~\eqref{ansatz_poly} 
corresponds to an inverse-power behavior of the 
remainder term,
\begin{equation}
\label{remainder_poly}
r_n = 
\sum_{j=1}^{\infty} c_j \, \left( \frac{1}{n+1} \right)^j
=  \frac{c_1}{n+1} + \frac{c_2}{(n+1)^2} + \dots \,.
\end{equation}
This, in turn, corresponds to Eq.~(7.1-2) of Ref.~\cite{We1989}.

%
%
\subsection{Neville Algorithm}
\label{sec2B}

For the calculation of the 
coefficient $c_0(n)$ of the 
interpolating polynomial, one
can use the Neville algorithm \cite{Ne1934}. 
Given $n+1$ partial sums $s_i$ of the input series
$\{ s_i \}_{n=0}^\infty$, the Neville 
algorithm refines the approximation of $s_\infty$ by 
iteratively interpolating between points defined 
by the sequence $s_i$. Recursive schemes have been 
described in Refs.~\cite{GaGu1974,vT1994}.
Adapted to our choice of abscissas $x_i = 1/(i+1)$, 
one may formulate a lozenge scheme as follows.
One starts the recursion with the values
\begin{equation}
s^0_i = s_i, \qquad i \in (0, \ldots, n) \,,
\label{eq:neville1}
\end{equation}
and uses a three-term recursion,
\begin{equation}
\label{three_rec}
s^m_i = \frac{(i+1) s^{m-1}_i - (i+1-m) \, s^{m-1}_{i-1}}{m},  
\qquad
1 \leq m \leq n \,,
\qquad
m \leq i \leq n \,.
\end{equation}
One increases $m$, thereby decreasing the allowable 
values of $i$ in the process, until,
for $m=n$, only one possible value is
left for $i$, namely, $i = m = n$. Finally,
the desired coefficient $c_0(n)$ is recovered
as follows,
\begin{equation}
c_0(n) = s^n_n \,.
\end{equation}
While the recursive scheme is very helpful, it does not establish a
direct connection between each term of the partial sums used for extrapolation
and the final estimated limit. Hence, no additional
information can be derived for the 
asymptotic behavior of the 
partial sums $s_n$ for large $n$.
In particular, while the recursive scheme enables the 
estimation of the limit $\lim_{x_i \to 0^+} f(x_i) = s_\infty$, 
it does not provide insights
into the behavior of the interpolating function $f(x_i)$
near $x_i=0$, e.g., the rate at which the
limit is approached.

Let us briefly clarify a notation used in the following. 
Our goal is to estimate $s_\infty$ utilizing the initial $n$
terms of the series. We denote the estimate of the 
limit after the convergence acceleration transformation 
by the symbol $\calT_n$.
Let now $\calT$ be the transformation 
which takes as input the partial sums $(s_0, s_1, \ldots, s_n)$,
and by which one obtains $\calT_n$,
\begin{equation}\label{eq:Hn}
\calT_n = \calT(s_0, s_1, \ldots, s_n).
\end{equation}
For example, if $\calT$ is the Neville algorithm 
as given in the adaptation~\eqref{three_rec},
then $\calT_n = c_0(n) = s_n^n$.
For the other methods discussed here, 
the definition of $\calT_n$ is adapted in the 
obvious way.

%
%
\subsection{Universal Formula}
\label{sec2C}

In order to obtain universal formulas which 
allow access to subleading terms (that serve to 
measure the rate of convergence of the input series),
it is advantageous to map the problem
onto a system of linear equations.
One interprets the 
coefficients $c_j(n)$ of the polynomial 
as the elements of a vector $\vec c(n) = ( c_0(n), ..., c_n(n) )$
with $n+1$ elements,
\begin{equation}
s_i = \sum_{j=0}^n c_j(n) \,
\left( \frac{1}{i + 1} \right)^j 
= \sum_{j=0}^n M_{ij}(n) \, c_j(n) \,,
\qquad
M_{ij}(n) = \left( \frac{1}{i+1} \right)^j \,,
\qquad i, j \in (0, \ldots, n) \,.
\end{equation}
The elements $M_{ij}(n)$ can be interpreted as components of a
matrix of dimension $(n+1) \times (n+1)$,
in which case one has $\vec s = \mathbbm{M}(n) \cdot \vec c(n)$, 
where $\vec s = ( s_0, ..., s_n )$ is the vector 
of the first $n+1$ partial sums of the input series.
Therefore, the solution can be obtained by matrix inversion,
\begin{equation}
\vec c(n) = [\mathbbm{M}(n)]^{-1} \cdot \vec s \,,
\qquad
\vec c = ( c_0(n), ..., c_n(n) ) \,.
\end{equation}
The inversion of the matrix $\mathbbm{M}(n)$
can lead to significant numerical instability and loss of
accuracy. This is because large numbers can cause overflow or underflow issues
in floating-point arithmetic, while the entries 
of $\mathbbm{M}(n)$ span many orders of magnitude,
which leads to unfavorable condition numbers for higher $n$.
This can exacerbate rounding errors. For this reason, it is 
advantageous to derive analytic formulas for 
the first five rows of the inverted matrix $[\mathbbm{M}(n)]^{-1}$,
which leads to analytic formulas 
for the coefficients $c_0(n)$, $c_1(n)$, $c_2(n)$, $c_3(n)$,
and $c_4(n)$. We find that 
\begin{subequations}
\label{turboneville}
\begin{align}
c_j(n) =& \; \sum_{i=0}^n \frac{(-1)^{n-i} \, (i+1)^n}%
{\Gamma(i+1) \, \Gamma(n-i+1)} \, \calP_j(i,n) \, s_i\,,
\qquad j \in (0, 1, 2, 3, 4) \,,
\\
\calP_0(i,n) =& \; 1 \,,
\\
2 \, \calP_1(i,n) =& \; 2 i - 3 n - n^2 \,,
\\ 
24 \calP_2(i,n) =& \; 
24 i + 24 i^2 - 26 n - 36 i n + 9 n^2 - 12 i n^2 + 14 n^3 + 3 n^4 \,,
\\
48 \calP_3(i,n) =& \; 48 i + 96 i^2 + 48 i^3 - 48 n - 
   124 i n - 72 i^2 n + 26 n^2 
\nonumber\\
& \; -  6 i n^2 - 24 i^2 n^2 + 29 n^3 + 
   28 i n^3 - n^4 + 6 i n^4 - 5 n^5 - n^6 \,.
\\
5760 \calP_4(i,n) =& \; 
  5760 i + 17280 i^2 + 17280 i^3 + 5760 i^4 - 5712 n - 
  20640 i n - 23520 i^2 n  - 8640 i^3 n 
\nonumber\\
& \; 
  + 3380 n^2 + 2400 i n^2 - 3600 i^2 n^2 - 2880 i^3 n^2 + 3660 n^3 + 6840 i n^3 
  + 3360 i^2 n^3 
\nonumber\\
& \; 
  - 385 n^4 + 600 i n^4 + 720 i^2 n^4 - 888 n^5 
- 600 i n^5 - 130 n^6 - 120 i n^6 + 60 n^7 + 15 n^8 \,.
\end{align}
\end{subequations}
Results for the coefficients $c_j(n)$ 
with $j > 4$ are given in Appendix~\ref{appa}.
The results given in Eq.~\eqref{turboneville}
facilitate the estimation of the 
functional relationship of the $f(x_i)$ for 
small abscissas $x_i$, thereby providing 
insight the convergence rate of the partial sums 
of the input series toward its limit.
Specifically, one notes the relationship
\begin{align}
s_n = & \; s_\infty + r_n \approx c_0 + \frac{c_1}{n + 1} + \frac{c_2}{(n + 1)^2} +
\dots \,, \qquad n \to \infty \,,
\qquad
c_j \equiv \lim_{n\to \infty} c_j(n) \,,
\\
r_n = & \; s_n - s_\infty \approx \frac{c_1}{n + 1} 
+ \frac{c_2}{(n + 1)^2} + \dots \,, \qquad n \to \infty \,,
\end{align}
where careful attention is given to the sign of the remainder term $r_n$.

\begin{figure}[t!]
\begin{center}
\begin{minipage}{0.8\linewidth}
\begin{center}
\includegraphics[width=0.7\linewidth]{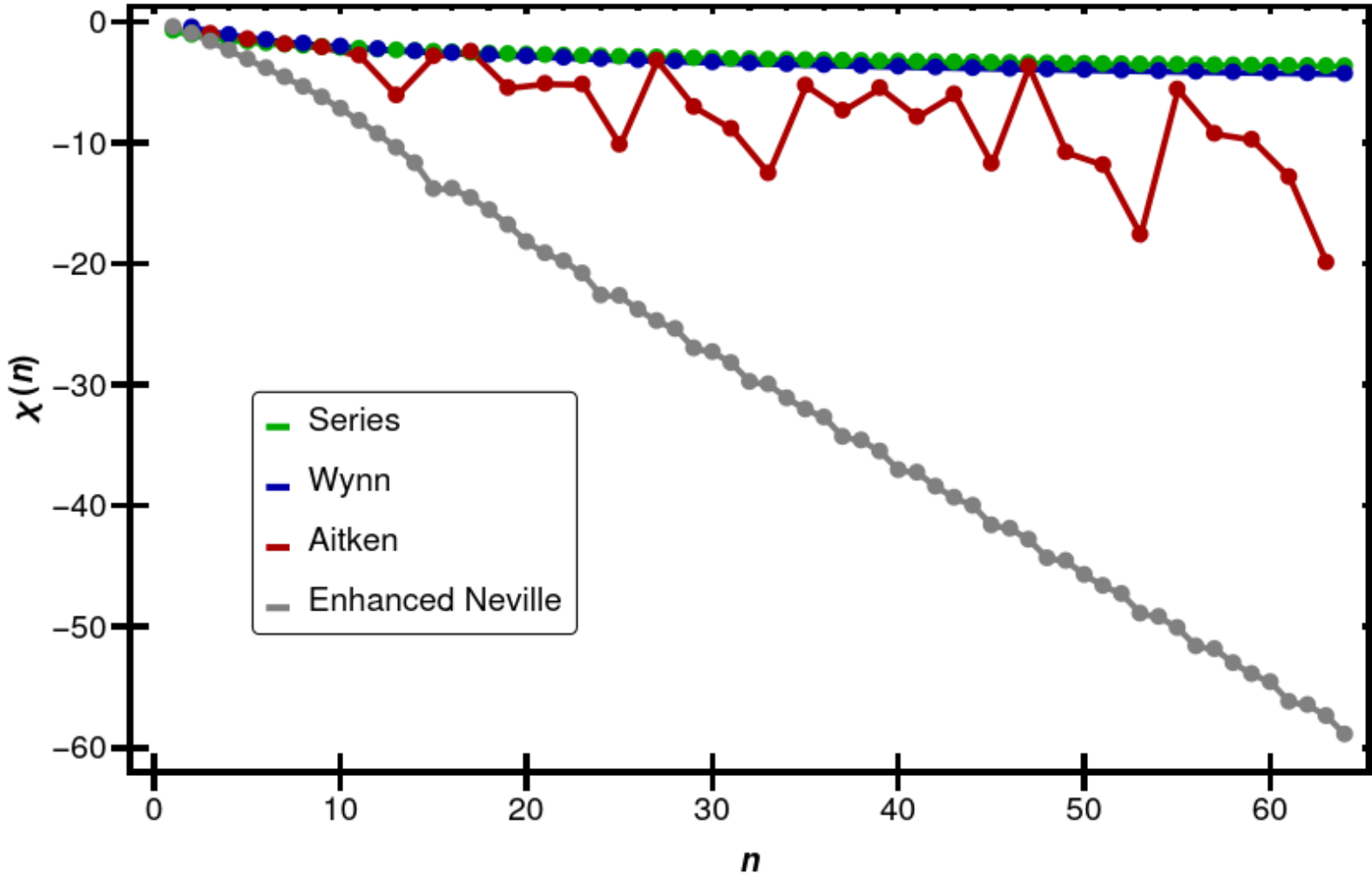}
\end{center}
\caption{\label{fig1}  
We investigate the convergence of the enhanced Neville transformation (green
curve) for the input series~\eqref{example}, as measured by the quantity
$\chi(n)$ defined in Eq.~\eqref{chi}. The quantity $\chi(n)$ roughly decreases
by unity with every iteration of the algorithm, indicating the gain of roughly
one converged decimal with every order of the 
enhanced Neville transformation. The
convergence of the enhanced Neville transformation is compared with the convergence
of the Wynn epsilon algorithm
(blue curve), of the Aitken $\Delta^2$ process
(red curve), and of the series 
itself (black curve).
(For interpretation of the colors in the ﬁgure, the reader is referred to the
current web version.
The curves labeled ``Series'', ``Wynn'', ``Aitken'' and ``Enhanced Neville''
are listed from top to bottom both in the figure legend as well as
in the figure itself.)}
\end{minipage}
\end{center}
\end{figure}

\begin{figure}[t!]
\begin{center}
\begin{minipage}{0.8\linewidth}
\begin{center}
\includegraphics[width=0.8\linewidth]{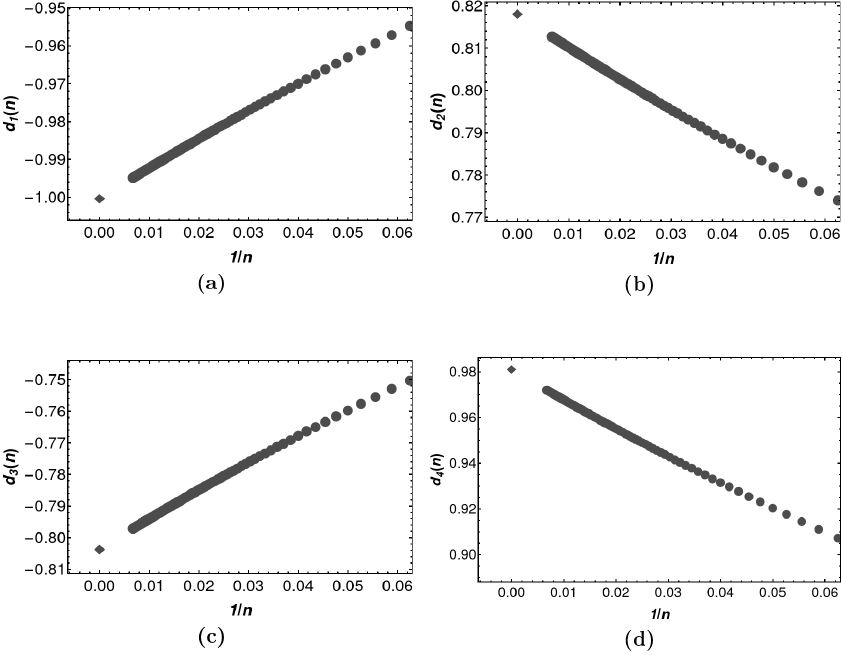}
\end{center}
\caption{\label{fig2} The figures 
illustrate the convergence of sequences 
$d_j(n)$ toward $c_j$ for $j = 1, 2, 3, 4$,
according to Eq.~\eqref{eq:sigmas}.
The cases $j = 1, 2, 3, 4$
are treated in Fig.~(a), (b), (d), and (d),
respectively.
The first 150 $d_j(n)$ with $n = 0, \dots, 149$
are compared to the value $c_j = d_j(\infty)$,
the latter being denoted by a black diamond.
One notes that $c_1 = -1$ is integer-valued,
while the first 50 decimals of the 
coefficients $c_2$, $c_3$ and $c_4$ 
are given in Eq.~\eqref{eq:coeffs}.}
\end{minipage}
\end{center}
\end{figure}

%
%
\section{Comparison to the Other Methods}
\label{sec3}

Some remarks may be in order with respect to the comparison of the method
outlined above to other established methods like the Aitken's $\Delta^2$
process and Wynn's epsilon algorithm~\cite{We1989,BrRZ1991,CaEtAl2007}.
These are two of the most well-known and frequently used
algorithms for series convergence acceleration. Both of them are based on the
Shanks transformation \cite{Sh1955}, which produces Padé approximants if the
input data are the partial sums of a power series~\cite{We1989}. The former produces an
efficient recursive scheme to compute conveniently such transformation. The
formula for the epsilon algorithm is given by:

\begin{equation}
\epsilon_{k+1}^{(n)} = \epsilon_{k-1}^{(n+1)} + \frac{1}{\epsilon_{k}^{(n+1)} - 
\epsilon_{k}^{(n)}}, 
\qquad \mbox{$-1 < k < n - 1$}  \,,
\end{equation}
with $k$ indicating the iteration level and $n$ the position within that level.
The initial conditions are $\epsilon_{-1}^{(n)} = 0$ for all $n$, and
$\epsilon_{0}^{(n)} = s_n$, where $s_n$ is the $n$th term of the sequence being
accelerated.
The values for $n=0$, $\epsilon^{(0)}_k$, describe a sequence
of upper-diagonal (odd $k$) and diagonal (even $k$) Pad\'{e} approximants to
the input series. The result of the transformation of the partial sums
$(s_0,\dots, s_n)$ is written using the diagonal Pad\'{e} approximants, that is
$\calT_n = \epsilon^{(0)}_{\floor*{\frac{n}{2}}}$, with $\floor*{\cdot}$ the
integer part.
Aitken's $\Delta^2$ process consists instead of a recursive
iteration of the first-order Shanks transformation,
\begin{equation}
A_{k+1}^{(n)} = A_{k}^{(n)} - 
\frac{\left(\Delta A_{k}^{(n)}\right)^2}{\Delta^2 A_{k}^{(n)}},
\end{equation}
where, as for the epsilon algorithm, 
$k$ indicates the iteration level and $n$ the 
position within that level. The initial conditions are 
$A_{0}^{(n)}=s_n$, while 
$\Delta A_{k}^{(n)} = A_{k}^{(n+1)} - A_{k}^{(n)}$ is the first
difference of the sequence, 
and $\Delta^2 A_{k}^{(n)} = \Delta A_{k}^{(n+1)} - \Delta A_{k}^{(n)}$ is
the second-order difference of the sequence.
The result of the transformation of the 
partial sums $(s_0,\dots, s_n)$ then is 
$\calT_n = A^{(0)}_{\floor*{\frac{n-1}{2}}}$, 
with $\floor*{\cdot}$ denoting the integer part.
A meaningful iteration of the $\Delta^2$ process,
which leads to an improved convergence,
therefore is attained when $n$ is increased by two.

Our assumptions, formulated in
Eq.~(\ref{remainder}), imply that the input series exhibits logarithmic
convergence.  Specifically, one has the asymptotic relation
\begin{equation}
\lim_{n\to\infty}\frac{s_{n+1}-s}{s_{n}-s}=\lim_{n\to\infty}\frac{r_{n+1}}{r_{n}}=1.
\end{equation}
This observation indicates that the Shanks transform, which is predicated on
the assumption that, in the sequence's tail (for large $n$), the ratio
\begin{equation}\label{eq:approx}
\frac{s_{n+1}-s}{s_{n}-s}\approx\frac{s_{n+2}-s}{s_{n+1}-s}\approx \lambda,
\end{equation}
remains constant, may not yield accurate results. 
In fact, advantage can be taken in Eq.~(\ref{eq:approx}) 
for a linearly convergent series (i.e., for the case 
$\lambda<1$). In this case, the sequence
terms are expected to follow
\begin{equation}
s_n=\sigma_0+\sigma_1 \alpha_1^n+O(\alpha_2^n), \quad |\alpha_2|<|\alpha_1|,
\end{equation}
indicating that the error in approximating Eq.~(\ref{eq:approx}) 
diminishes geometrically with $n$. Specifically,
\begin{equation}
\frac{s_{n+1}-s}{s_{n}-s}-\frac{s_{n+2}-s}{s_{n+1}-s} =
\calO\left(\frac{\alpha_2}{\alpha_1}\right)^n.
\end{equation}
However, for a series whose terms are represented by Eq.~(\ref{ansatz_poly}), the
error reduction is merely polynomial in $n$:
\begin{equation}
\frac{s_{n+1}-s}{s_{n}-s}-\frac{s_{n+2}-s}{s_{n+1}-s}=\calO\left(\frac{1}{n}\right),
\end{equation}
which leads to a much slower rate 
of convergence for large $n$ as compared to 
the factor $\lambda = \alpha_2 / \alpha_1 < 1$.
Consequently, methodologies dependent on 
such algorithms as Aitken's $\Delta^2$ process~\cite{Ai1926},
the Shanks transformation~\cite{Sh1955}, 
or Wynn's $\epsilon$ algorithm~\cite{Wy1956eps},
demonstrate suboptimal performance for these series.
For our particular model example,
given in Eq.~\eqref{example},
these statements are illustrated in Fig.~\ref{fig1}.
A clarifying remark is in order.
For the comparison in Fig.~\ref{fig1}, we plot the 
data points as a function of $n$, where $n$ is the 
maximum index of the terms of the input series used
in order to calculate the convergence acceleration
transform. We indicate the error of the extrapolation as a function
of $n$. Some algorithms, such as the Aitken method,
only produce a sensible result when
$n$ is being increased in steps of two, $n \to n+2$.
This explains why there are some apparent ``gaps'' in the
numerical data.

%
%
\section{Numerical Examples}
\label{sec4}

%
%
\subsection{Slowly Convergent Model Series}
\label{sec4A}

We discuss the application of the enhanced Neville 
algorithm to the infinite series with partial sums
\begin{equation}
\label{example}
s_n = \sum_{k=0}^n a_k \,,
\qquad
a_k = \frac{4}{\pi \, (k+1)^2} \,
\arctan\left( \frac{k+2}{k+3} \right) \,,
\end{equation}
where the prefactor is chosen so that 
$a_k \approx 1/k^2$ for $k \to \infty$.
We calculate the enhanced Neville transformations, 
up to order $n=100$, and the determination of the convergence 
according to the formula
\begin{equation}
\label{chi}
\chi(n) = 
\frac{ \ln( | \calT_n - \calT_{n + 1} | ) }{\ln(10)} \,,
\end{equation}
with $\calT_n$ defined in Eq. (\ref{eq:Hn}).
The quantity $\chi(n)$ measures the apparent convergence 
of the transforms as the order of the 
transformation is increased, by calculating the 
number of apparently converged decimal digits.

In Fig.~\ref{fig1}, we compare the convergence of the 
enhanced Neville transformation
to Aitken's $\Delta^2$ process~\cite{Ai1926},
Wynn's epsilon algorithm~\cite{Wy1956eps} and 
alongside the series itself. The enhanced Neville scheme
is shown to outperform the other convergence acceleration algorithms, securing
approximately one additional decimal of convergence with each transformation
order. In contrast, the other methods progress much more slowly.
This observation is tied to the asymptotic structure of the 
input series, according to Eq.~\eqref{asymptotic}.

Our analysis, based on the enhanced Neville algorithm,
yields a 50-decimal approximation of the series limit:
\begin{equation}
s_\infty = \lim_{n \to \infty} c_0(n) = 
1.31279\,49538\,25865\,79196\,34865\,83906\,40442\,73891\,27574\,77554.
\end{equation}
Furthermore, we obtain the results
\begin{subequations}
\label{eq:coeffs}
\begin{align}
c_1 = \lim_{n \to \infty} c_1(n) =& \; -1  \,, \qquad
c_2 = \lim_{n \to \infty} c_2(n) = \frac12 + \frac{1}{\pi} =
 0.81830\,98861\,83790\dots\,,
\\
c_3 = \lim_{n \to \infty} c_3(n) =& \; -\frac16 - \frac{2}{\pi} = 
-0.80328\,64390\,34248\dots \,,
\qquad
c_4 = \lim_{n \to \infty} c_4(n) = \frac{37}{12 \pi} = 
 0.98145\,54824\,00021\dots\,.
\end{align}
\end{subequations}
The next higher coefficients
are $c_5 = \frac{1}{30} - \frac{131}{30 \pi}$,
$c_6 = \frac{268}{45 \pi}$,
$c_7 = -\frac{1}{42} - \frac{549}{70 \pi}$,
$c_8 = \frac{98347}{10000 \pi}$,
$c_9 = -\frac{1}{30} - \frac{9253}{840 \pi}$,
and $c_{10} = \frac{66011}{6300 \pi}$.
These results have been verified to about 100 digits
against numerical data obtained with the help of the 
enhanced Neville algorithm described here. They can be 
obtained analytically based on the formula 
\begin{equation}
s_n = \sum_{k=0}^n a_k = s_\infty - \sum_{k=n+1}^\infty a_k  \,,
\end{equation}
while realizing that the terms $a_k$ in the latter sum from 
$k=n+1$ to $k = \infty$ can be expanded for large summation index $k$.
We define the approximations $d_{j}(n) \approx c_j$ as follows,
\begin{equation}
\label{eq:sigmas}
d_{j}(n) = (n+1)^j 
\left( s_n - \sum_{r=0}^{j-1}\frac{c_r}{(n+1)^r} \right), 
\qquad 
\lim_{n \to \infty} d_{j}(n) = c_j \,,
\qquad j = 1,2,3,4\,.
\end{equation}
It is instructive to plot the sequence
of the $d_{j}(n)$ against $1/n$, for $j=1,2,3,4$,
as shown in Fig.~\ref{fig2},
for the model series given in Eq.~\eqref{example}.
The value at infinite $n$ (where $1/n \to 0$) is
given by $\lim_{n \to \infty} d_{j}(n) = c_j$.
The trend of the data presented in Fig.~\ref{fig2}
is consistent with the convergence of the approximants 
$d_j(n) \to c_j$ for $j=1,2,3,4$, for large $n$.

%
%
\subsection{Bethe Logarithm}
\label{sec4B}

We investigate the series of partial sums $s_n$, where
\begin{equation}
\label{blog}
b_n = \sum_{k=0}^n a_k \,,
\qquad
a_k = \frac{16(k+2)}{(k+1)^2 \, (k+3)^2} \,
\Phi\left( -\frac{3 + k}{1 + k}, 1, 2k+4 \right)  \,,
\end{equation}
where $\Phi$ is the Lerch phi transcendent,
with $\Phi(z,s,a) = \sum_{n=0}^\infty z^n/(n+a)^s$.
According to Chap.~4 of Ref.~\cite{JeAd2022book},
the limit $s_\infty$ of the input series 
is related to the Bethe logarithm~\cite{Be1947} 
$\ln k_0(1S)$ of the hydrogen ground-state
by the following relation,
\begin{equation}
\label{lnk0}
b_\infty = \ln k_0(1S) - 10 \, \ln(2) + 2 \, \zeta(2) + 1 \,.
\end{equation}
The sum $b_\infty$ and hence, the Bethe logarithm 
$\ln k_0(1S)$ of the ground state of hydrogen,
is calculated in the same manner as the 
example treated in Sec.~\ref{sec4A}.
The enhanced Neville algorithm obtains the 50-figure result
\begin{equation}
\label{lnk01S}
\ln k_0(1S) = 2.98412\,85557\,65497\,61075\,97770\,90013\,79796\,99751\,80566\,17002
\end{equation}
in the 58th order of transformation. It is somewhat surprising that
the coefficients $c_j$ with $j=1,2,3,4$ are rational numbers,
\begin{equation}\begin{split}
\label{miracle}
c_1=&\,0, \qquad c_2=\,0, \qquad c_3=-\frac{4}{3}, \qquad c_4=\frac{27}{4}, 
\qquad c_5= -\frac{703}{30}, \qquad c_6=\frac{3329}{48}, \\
c_7=& -\frac{63163}{336}, \qquad c_8=\frac{184961}{384}, \qquad
c_9=-\frac{569323}{480}, \qquad c_{10}=\frac{7256477}{2560} \,.
\end{split}\end{equation}
The results given in Eq.~\eqref{miracle} currently constitute conjectures which
we have verified to at least 100 decimals. The result
$c_1 = c_2 = 0$ is expected because $\lim_{k \to \infty}a_k
\sim \frac{1}{k^4}$, which then leads to $\lim_{k \to \infty}b_k = b_{\infty} +
\mathcal{O}\left(\frac{1}{k^3}\right)$. Values of the Bethe
logarithms $\ln k_0(1S)$, $\ln k_0(2S)$ and $\ln k_0(2P)$ accurate to 100
digits are presented in Appendix~\ref{appb}.  These could be used for the
search of possible analytic formulas using the PSLQ
algorithm~\cite{FeBa1992,BaPl1997,FeBaAr1999,BaBr2001}.  The PSLQ algorithm can
be used to search for analytic expressions of accurately known numerical
quantities in terms of known constants such as the Euler constant $\gamma_E =
0.57721\dots$, various Riemann zeta functions, powers of $\pi$, and
multiplicative combinations of these constants.  We can report that we have
carried out a limited set of searches with the same constants that were used in
Eq.~(A11) of Ref.~\cite{GiEtAl2020} without success.

%
%
\section{Generalizations of the Algorithm}
\label{sec5}

%
%
\subsection{Generalizations to Half--Integer--Power Remainder Terms}
\label{sec5A}

It is instructive to briefly consider possible 
generalizations. 
For example, we may easily adapt the algorithm to series
with half-integer powers, which leads to the {\em ansatz}
\begin{equation}
f(x_i) = \sum_{j=0}^{\infty} c^\half_j \, x^\half_{ij}
\approx \sum_{j=0}^{n} c^\half_j \, x^\half_{ij}
= c^\half_0 + \frac{c^\half_1}{\sqrt{i+1}} + \frac{c^\half_2}{i+1} +
\ldots + \frac{c^\half_n}{(i+1)^{n/2}} \,,
\label{eq:polyHI}
\end{equation}
where the $x^\half_{ij}$ are given as follows,
\begin{equation}
\label{absHI}
x^\half_{ij} = \frac{1}{(i+1)^{j/2}} \,.
\end{equation}
Our modified remainder term
is thus based on half-integer powers,
\begin{equation}
r^\half_n = \frac{c^\half_1}{(n+1)^{1/2}} 
+ \frac{c^\half_2}{(n+1)} + \frac{c^\half_3}{(n+1)^{3/2}} + \dots \,.
\end{equation}
Given $n+1$ partial sums $s_i$ of the input series
$\{ s_i \}_{n=0}^\infty$, the generalized algorithm
aims to solve the system of equations
\begin{equation}
\label{linear}
s_i = \sum_{j=0}^n c^\half_j(n) \, x^\half_{ij}
\qquad i, j \in (0, \ldots, n) \,,
\end{equation}
as a function of $n$, leading to 
results for the coefficients $c^\half_j(n)$ with 
$j \in (0, \ldots, n)$. In full analogy to 
Eqs.~\eqref{conv_poly} and~\eqref{conv_factorial},
one then assumes that the
value of the infinite series is recovered as
$n$ is increased,
\begin{equation}
\label{conv_half_integer}
s = s_\infty = \lim_{n \to \infty} c^\half_0(n) \,.
\end{equation}
Using the half-integer-power Neville scheme, we investigate the 
acceleration of the convergence of the following series,
\begin{equation}
\label{example_sqrt}
s_n = \sum_{k=0}^n a_k \,,
\qquad
a_k = \frac{4}{\pi \, (\sqrt{k}+1)^4} \,
\arctan\left( \frac{\sqrt{k} + 2}{\sqrt{k} + 3} \right) \,,
\end{equation}
where the prefactor is adjusted so 
that $a_k \approx 1/k^2$ for $k \to \infty$.
The transforms $c_0(n \leq 203)$ converge to the 
100-decimal result
\begin{align}
s_\infty = \lim_{n \to \infty} c^\half_0(n) = 
0.&92316\,76494\,43269\,09201\,53242\,74308\,14908\,92138\,92222\,83308\nonumber\\
  &78498\,38145\,30824\,82489\,27413\,45388\,44628\,47523\,62540\,58525 \,.
\end{align}
The transforms are calculated by solving the linear system 
given in Eq.~\eqref{linear} by standard matrix 
inversion techniques~\cite{PrFlTeVe1993}.
Furthermore, the following results are obtained,
\begin{equation}
c^\half_1 = 0\,, \qquad c^\half_2 = -1 \,, \qquad
c^\half_3 = \frac13 \left( 8 + \frac{4}{\pi} \right) \,.
\end{equation}
These are exact.
By contrast, the convergence acceleration method outlined 
in Sec.~\ref{sec2} has only converged to about nine decimals
for $n=200$.

%
%
\subsection{Generalizations to Fractional--Power Remainder Terms}
\label{sec5B}

In Sec.~\ref{sec5B}, we had considered the 
generalization of the enhanced Neville
algorithm to series whose remainder terms 
comprise half-integer powers. 
Here, shall generalize the algorithm
further, to cases where the remainder
terms can be expressed as powers
of $n^{-1/s}$. We shall provide an example
for the case $s=3$. The {\em ansatz} is
\begin{equation}
f(x_i) = \sum_{j=0}^{\infty} c^\expo_j \, x^\expo_{ij}
\approx \sum_{j=0}^{n} c^\expo_j \, x^\expo_{ij}
= c^\expo_0 
+ \frac{c^\expo_1}{(i+1)^{1/s}} 
+ \frac{c^\expo_2}{(i+1)^{2/s}} +
\ldots + \frac{c^\expo_n}{(i+1)^{n/s}} \,,
\label{eq:polyFRAC}
\end{equation}
where the $x^\expo_{ij}$ are given as follows,
\begin{equation}
\label{absFRAC}
x^\expo_{ij} = \frac{1}{(i+1)^{j/2}} \,.
\end{equation}
The generalized {\em ansatz} for the remainder term is 
\begin{equation}
r^\expo_n = \frac{c^\expo_1}{(n+1)^{1/s}} 
+ \frac{c^\expo_2}{(n+1)^{2/s}} + \frac{c^\expo_3}{(n+1)^{3/s}} + \dots \,.
\end{equation}
Based on $n+1$ partial sums $s_i$ of the input series
$\{ s_i \}_{n=0}^\infty$, we aim to solve the 
system of $(n+1)$ equations
\begin{equation}
\label{linear_frac}
s_i = \sum_{j=0}^n c^\expo_j(n) \, x^\expo_{ij}
\qquad i, j \in (0, \ldots, n) \,,
\end{equation}
for the coefficients $c^\expo_j(n)$ with $j \in (0, \ldots, n)$. 
In full analogy to 
Eqs.~\eqref{conv_poly} and~\eqref{conv_factorial},
one then assumes that the
value of the infinite series is recovered as
$n$ is increased,
\begin{equation}
\label{conv_frac}
s = s_\infty = \lim_{n \to \infty} c^\expo_0(n) \,.
\end{equation}
Using the fractional-power Neville scheme
with $s=3$, we investigate the 
acceleration of the convergence of the following series,
\begin{equation}
\label{example_frac}
s_n = \sum_{k=0}^n a_k \,,
\qquad
a_k = \frac{4}{\pi \, (k^{1/3}+1)^6} \,
\arctan\left( \frac{k^{1/3} + 2}{k^{1/3} + 3} \right) \,,
\end{equation}
where the prefactor is adjusted so 
that $a_k \approx 1/k^2$ for $k \to \infty$.
The transforms $c^{(1/3)}_0(n \leq 402)$ converge to the result
\begin{align}
s_\infty = \lim_{n \to \infty} c^{(1/3)}_0(n) = 
0.&79903\,08857\,54877\,20184\,41296\,61307\,30047\,47629\,74541\,71193\nonumber\\
  &85133\,42481\,37440\,38454\,19813\,14657\,35529\,94968\,85946\,86410\,,
\end{align}
which has 100 decimals.
The transforms are calculated by solving the linear system 
given in Eq.~\eqref{linear_frac} by standard matrix 
inversion techniques~\cite{PrFlTeVe1993}.
Furthermore, the following results are obtained,
\begin{equation}
c^{(1/3)}_1 = c^{(1/3)}_2 = 0 \,, \qquad
c^{(1/3)}_3 = -1 \,.
\end{equation}
These are exact.
By contrast, the convergence acceleration method outlined 
in Sec.~\ref{sec2} has only converged to about ten decimals
for $n=400$.

%
%
\section{Conclusions}
\label{sec6}

For slowly convergent series which possess
an asymptotic expansion compatible with Eq.~\eqref{asymptotic},
the enhanced Neville algorithm described here
harvests the asymptotic structure to eliminate the 
maximum number of asymptotic terms in the 
expansion in inverse powers of $n$, 
in any given order to the transformation.
We take advantage of the asymptotic structure
of the terms to be summed, according to Eq.~\eqref{asymptotic}.
The general structure of the transformation
has been described in Sec.~\ref{sec2A} and
and compared to the recursive Neville algorithm 
in Sec.~\ref{sec2B}.
Our matrix-based formulation of the problem 
in Sec.~\ref{sec2A} enables us to not 
only derive formulas for the limit of the 
series, but also for subleading terms, 
according to Eq.~\eqref{turboneville}.

In Sec.~\ref{sec3},
we have compared the performance of the 
enhanced Neville algorithm to that 
of Wynn's $\epsilon$ algorithm~\cite{Wy1956eps},
which calculates Pad\'{e} approximants,
and Aitken's $\Delta^2$ process~\cite{Ai1926},
including the iterated 
Aitken process. We find that our method 
provides numerically superior results,
and harvests the asymptotic structure
of the input series, according to Eq.~\eqref{asymptotic}.
Furthermore, in comparison to the combined nonlinear-condensation
transformation, described in
Refs~\cite{JeMoSoWe1999,AkSaJeBeSoMo2003}, the enhanced Neville transformation
eliminates the need to ``sample''
the input series at very high orders, 
which is otherwise necessitated by the van Wijngaarden 
transformation~\cite{NP1961},
which is the first step of the combined nonlinear-condensation
transformation.

The versatility of the algorithm
is demonstrated by considering
generalizations to remainder terms 
involving inverse half-integer powers
(Sec.~\ref{sec5A}), general inverse fractional
powers (Sec.~\ref{sec5B}),
as well as remainder terms constituting inverse
factorial series (Appendix~\ref{secB1}),
as well as remainder terms containing only
even and quartic powers (see Appendix~\ref{secB2}).

For the Bethe logarithm, discussed in Sec.~\ref{sec4B},
we find that the subleading terms of the partial sums 
$b_n$ of the series given in Eq.~\eqref{blog} 
can be expressed in terms of rational coefficients,
\begin{equation}
b_n = \ln k_0(1S) - 10 \, \ln(2) + 2 \, \zeta(2) + 1 
- \frac{4}{3 (n+1)^3} - \frac{27}{4 (n+1)^4} + \calO\left(\frac{1}{n^5}\right) \,,
\end{equation}
for large $n$. For the definition of $b_n$, see Eq.~\eqref{blog}.
The coefficients of the subleading terms 
($-4/3$ and $-27/4$) have been determined 
based on our enhanced Neville algorithm.

Coefficients for the derivation of 
higher-order subleading asymptotics 
are given in Appendix~\ref{appa}.
Furthermore, Bethe logarithms for 
the $2S$ and $2P$ states are discussed in 
Appendix~\ref{appb}, and numerical results
are given to 100~decimal digits. 
With the enhanced Neville algorithm presented here, 
the accuracy of the calculation of the Bethe 
logarithms can easily be enhanced beyond the
100~decimals indicated in Appendix~\ref{appb},
facilitating a possible search for 
analytic representations using the PSLQ 
algorithm~\cite{FeBa1992,BaPl1997,FeBaAr1999,BaBr2001}.

\section{Acknowledgements}

This article is inspired by numerous insightful discussions
with Ernst Joachim Weniger on related problems,
dating back to the years 1999--2012.
The authors also acknowledge helpful conversations with Jean Zinn--Justin.
UDJ was supported by the National Science Foundation
(grant NSF PHY--2110294).
LTG gratefully acknowledges support from the Swedish Research Council
(Vetenskapsradet) Grant No. 638-2013-9243.

\appendix

\small

%
%
\section{Coefficients for Higher--Order Terms}
\label{appa}

We report in this Appendix the expressions for the 
polynomials $\calP_k(i,n)$ with $5 \leq k \leq 10$.
For $k=5$, the result is
\begin{align}
11520\calP_5(i,n) = & \; 11520 i + 46080 i^2 + 69120 i^3 + 46080 i^4 + 11520 i^5 
- 11520 n - 52704 i n - 88320 i^2 n \nonumber\\
& - 64320 i^3 n - 17280 i^4 n + 6768 n^2 + 11560 i n^2 - 2400 i^2 n^2 - 12960 i^3 n^2 - 5760 i^4 n^2 \nonumber\\
& + 7652 n^3 + 21000 i n^3 + 20400 i^2 n^3 + 6720 i^3 n^3 
- 680 n^4 + 430 i n^4 + 2640 i^2 n^4 + 1440 i^3 n^4 \nonumber\\
& - 2085 n^5 - 2976 i n^5 - 1200 i^2 n^5 - 395 n^6 - 500 i n^6 - 
240 i^2 n^6 + 198 n^7 + 120 i n^7 \nonumber\\
& + 70 n^8 + 30 i n^8 - 5 n^9 - 3 n^{10} \,.
\end{align}
For $k=6$, we have the result 
\begin{align}
2903040\calP_6(i,n) = & \; 
2903040 i + 14515200 i^2 + 29030400 i^3 + 29030400 i^4 + 14515200 i^5 + 2903040 i^6 \nonumber\\
& - 2914560 n - 16184448 i n - 35538048 i^2 n - 38465280 i^3 n - 20563200 i^4 n - 4354560 i^5 n \nonumber\\
& + 1667232 n^2 + 4618656 i n^2 + 2308320 i^2 n^2 - 3870720 i^3 n^2 - 4717440 i^4 n^2 - 1451520 i^5 n^2 \nonumber\\
& + 1942136 n^3 + 7220304 i n^3 + 10432800 i^2 n^3 + 6834240 i^3 n^3 + 1693440 i^4 n^3 - 97020 n^4 \nonumber\\
& - 63000 i n^4 + 773640 i^2 n^4 + 1028160 i^3 n^4 + 362880 i^4 n^4 - 523446 n^5 - 1275372 i n^5 \nonumber\\
& - 1052352 i^2 n^5 - 302400 i^3 n^5 - 146727 n^6 - 225540 i n^6 - 186480 i^2 n^6 - 60480 i^3 n^6 \nonumber\\
& + 44070 n^7 + 80136 i n^7 + 30240 i^2 n^7 + 30177 n^8 + 25200 i n^8 + 7560 i^2 n^8 + 406 n^9 \nonumber\\
& - 1260 i n^9 - 2205 n^{10} - 756 i n^{10} - 126 n^{11} + 63 n^{12} \,.
\end{align}
For $k=7$, the result is
\begin{align}
5806080\calP_7(i,n) = & \; 
5806080 i + 34836480 i^2 + 87091200 i^3 + 116121600 i^4 + 87091200 i^5 + 34836480 i^6 \nonumber\\
& + 5806080 i^7 - 5806080 n - 38198016 i n - 103444992 i^2 n - 148006656 i^3 n - 118056960 i^4 n \nonumber\\
& - 49835520 i^5 n - 8709120 i^6 n + 3303936 n^2 + 12571776 i n^2 + 13853952 i^2 n^2 - 3124800 i^3 n^2 \nonumber\\
& - 17176320 i^4 n^2 - 12337920 i^5 n^2 - 2903040 i^6 n^2 + 3783456 n^3 + 18324880 i n^3 \nonumber\\
& + 35306208 i^2 n^3 + 34534080 i^3 n^3 + 17055360 i^4 n^3 + 3386880 i^5 n^3 - 147000 n^4 \nonumber\\
& - 320040 i n^4 + 1421280 i^2 n^4 + 3603600 i^3 n^4 + 2782080 i^4 n^4 + 725760 i^5 n^4 - 920780 n^5 \nonumber\\
& - 3597636 i n^5 - 4655448 i^2 n^5 - 2709504 i^3 n^5 - 604800 i^4 n^5 - 313950 n^6 - 744534 i n^6 \nonumber\\
& - 824040 i^2 n^6 - 493920 i^3 n^6 - 120960 i^4 n^6 + 27573 n^7 + 248412 i n^7 + 220752 i^2 n^7 \nonumber\\
& + 60480 i^3 n^7 + 65187 n^8 + 110754 i n^8 + 65520 i^2 n^8 + 15120 i^3 n^8 + 14457 n^9 - 1708 i n^9 \nonumber\\
& - 2520 i^2 n^9 - 5397 n^{10} - 5922 i n^{10} - 1512 i^2 n^{10} - 1729 n^{11} - 252 i n^{11} 
+ 273 n^{12} + 126 i n^{12} + 63 n^{13} - 9 n^{14}.
\end{align}
With $C_8 = 1393459200$, we have the following formula for $k=8$,
\begin{align}
C_8 \, \calP_8(i,n) = & \; 1393459200 i + 9754214400 i^2 + 29262643200 i^3 + 48771072000 i^4 
+ 48771072000 i^5 + 29262643200 i^6 \nonumber\\
& + 9754214400 i^7 + 1393459200 i^8 - 1387653120 n 
- 10560983040 i n - 33994321920 i^2 n \nonumber\\
& - 60348395520 i^3 n - 63855267840 i^4 n - 40294195200 i^5 n - 14050713600 i^6 n - 2090188800 i^7 n \nonumber\\
& + 807277824 n^2 + 3810170880 i n^2 + 6342174720 i^2 n^2 + 2574996480 i^3 n^2
- 4872268800 i^4 n^2 \nonumber\\
& - 7083417600 i^5 n^2 - 3657830400 i^6 n^2 - 696729600 i^7 n^2 
+ 891826560 n^3 + 5306000640 i n^3 \nonumber\\
& + 12871461120 i^2 n^3 + 16761669120 i^3 n^3 + 12381465600 i^4 n^3 + 
4906137600 i^5 n^3 + 812851200 i^6 n^3 \nonumber\\
& - 67571600 n^4 - 112089600 i n^4 + 264297600 i^2 n^4 + 1205971200 i^3 n^4 + 1532563200 i^4 n^4 \nonumber\\
& + 841881600 i^5 n^4 + 174182400 i^6 n^4 - 204569280 n^5 - 1084419840 i n^5 - 1980740160 i^2 n^5 \nonumber\\
& - 1767588480 i^3 n^5 - 795432960 i^4 n^5 - 145152000 i^5 n^5 
 49594888 n^6 - 254036160 i n^6 - 376457760 i^2 n^6 \nonumber\\
& - 316310400 i^3 n^6 - 147571200 i^4 n^6 - 29030400 i^5 n^6 
- 1304520 n^7 + 66236400 i n^7 + 112599360 i^2 n^7 \nonumber\\
& + 67495680 i^3 n^7 + 14515200 i^4 n^7 + 6310455 n^8 + 42225840 i n^8 + 42305760 i^2 n^8 
+ 19353600 i^3 n^8 \nonumber\\
& + 3628800 i^4 n^8 + 5741280 n^9 + 3059760 i n^9 - 1014720 i^2 n^9 - 604800 i^3 n^9
+ 383204 n^{10} \nonumber\\
& - 2716560 i n^{10} - 1784160 i^2 n^{10} - 362880 i^3 n^{10} - 825840 n^{11} - 475440 i n^{11} 
- 60480 i^2 n^{11} - 76790 n^{12} \nonumber\\
& + 95760 i n^{12} + 30240 i^2 n^{12} + 57120 n^{13} + 15120 i n^{13} + 1260 n^{14} - 2160 i n^{14} 
- 1800 n^{15} + 135 n^{16}.
\end{align}
With $C_9 = 2786918400$
and $\calP_9\equiv \calP_9(i,n)$,
the result can be written, for $k=9$, as the 
sum of two terms,
\begin{align}
C_9 \, \calP_9 = & \calT_1 + \calT_2 \,, 
\\
\calT_1 =& \;
\; 2786918400 i + 22295347200 i^2 + 78033715200 i^3 + 156067430400 i^4 
+ 195084288000 i^5 \nonumber\\
& + 156067430400 i^6 + 78033715200 i^7 + 22295347200 i^8 + 2786918400 i^9 - 2786918400 n \nonumber\\
& - 23897272320 i n - 89110609920 i^2 n - 188685434880 i^3 n - 248407326720 i^4 n 
- 208298926080 i^5 n \nonumber\\
& - 108689817600 i^6 n - 32281804800 i^7 n - 4180377600 i^8 n 
+ 1642567680 n^2 + 9234897408 i n^2 \nonumber\\
& + 20304691200 i^2 n^2 + 17834342400 i^3 n^2 - 4594544640 i^4 n^2 - 23911372800 i^5 n^2 
- 21482496000 i^6 n^2 \nonumber\\
& - 8709120000 i^7 n^2 - 1393459200 i^8 n^2 + 1839755520 n^3 
+ 12395654400 i n^3 + 36354923520 i^2 n^3 \nonumber\\
& + 59266260480 i^3 n^3 + 58286269440 i^4 n^3 + 34575206400 i^5 n^3 + 11437977600 i^6 n^3 
+ 1625702400 i^7 n^3 \nonumber\\
& - 193851264 n^4 - 359322400 i n^4 + 304416000 i^2 n^4 + 2940537600 i^3 n^4 
+ 5477068800 i^4 n^4 \nonumber\\
& + 4748889600 i^5 n^4 + 2032128000 i^6 n^4 + 348364800 i^7 n^4 
- 486633040 n^5 - 2577978240 i n^5 \nonumber\\
& - 6130320000 i^2 n^5 - 7496657280 i^3 n^5 - 5126042880 i^4 n^5 
- 1881169920 i^5 n^5 - 290304000 i^6 n^5 \,.
\end{align}
The second term $\calT_2$ is given as follows,
\begin{align}
\calT_2 =& \;
- 53606320 n^6 - 607262096 i n^6 - 1260987840 i^2 n^6 - 1385536320 i^3 n^6 - 927763200 i^4 n^6 \nonumber\\
& - 353203200 i^5 n^6 - 58060800 i^6 n^6 + 40638488 n^7 + 129863760 i n^7 + 357671520 i^2 n^7 \nonumber\\
& + 360190080 i^3 n^7 + 164021760 i^4 n^7 + 29030400 i^5 n^7 
- 6368192 n^8 + 97072590 i n^8 \nonumber\\
& + 169063200 i^2 n^8 + 123318720 i^3 n^8 + 45964800 i^4 n^8 + 7257600 i^5 n^8 
- 252565 n^9 + 17602080 i n^9 \nonumber\\
& + 4090080 i^2 n^9 - 3239040 i^3 n^9 - 1209600 i^4 n^9 
+ 5530785 n^{10} - 4666712 i n^{10} - 9001440 i^2 n^{10} \nonumber\\
& - 4294080 i^3 n^{10} - 725760 i^4 n^{10} - 77996 n^{11} - 2602560 i n^{11} - 1071840 i^2 n^{11} 
- 120960 i^3 n^{11} \nonumber\\
&- 873524 n^{12} + 37940 i n^{12} + 252000 i^2 n^{12} + 60480 i^3 n^{12} 
+ 34690 n^{13} + 144480 i n^{13} + 30240 i^2 n^{13} \nonumber\\
& + 61670 n^{14} - 1800 i n^{14} - 4320 i^2 n^{14} - 6212 n^{15} - 3600 i n^{15} 
- 1620 n^{16} + 270 i n^{16} + 315 n^{17} - 15 n^{18}.
\end{align}
Finally, with $C_{10} = 367873228800$
and $\calP_{10}\equiv \calP_{10}(i,n)$,
we have the following result for $k=10$,
\begin{align}
C_{10} \, \calP_{10} = & \; 
367873228800 i + 3310859059200 i^2 + 13243436236800 i^3 + 30901351219200 i^4 + 46352026828800 i^5 \nonumber\\
& + 46352026828800 i^6 + 30901351219200 i^7 + 13243436236800 i^8 
+ 3310859059200 i^9 + 367873228800 i^{10} \nonumber\\
& - 370660147200 n - 3522313175040 i n - 14917040455680 i^2 n - 36669077913600 i^3 n \nonumber\\
& - 57696244531200 i^4 n - 60285225369600 i^5 n - 41842514165760 i^6 n - 18608254156800 i^7 n \nonumber\\
& - 4813008076800 i^8 n - 551809843200 i^9 n 
+ 211314216960 n^2 + 1435825391616 i n^2 + 3899225696256 i^2 n^2 \nonumber\\
& + 5034352435200 i^3 n^2 
+ 1747653304320 i^4 n^2 - 3762781102080 i^5 n^2 - 5991990681600 i^6 n^2 \nonumber\\
& - 3985293312000 i^7 n^2 - 1333540454400 i^8 n^2 - 183936614400 i^9 n^2 
+ 252909144576 n^3 \nonumber\\
& + 1879074109440 i n^3 + 6435076285440 i^2 n^3 + 12621996288000 i^3 n^3 
+ 15516933949440 i^4 n^3 \nonumber\\
& + 12257714810880 i^5 n^3 + 6073740288000 i^6 n^3 + 1724405760000 i^7 n^3 
+ 214592716800 i^8 n^3 \nonumber\\
& 
- 12322523136 n^4 - 73018923648 i n^4 - 7247644800 i^2 n^4 + 428333875200 i^3 n^4 
+ 1111124044800 i^4 n^4 \nonumber\\
& + 1349826508800 i^5 n^4 + 895094323200 i^6 n^4 + 314225049600 i^7 n^4 + 45984153600 i^8 n^4 
- 77131325216 n^5 \nonumber\\
& - 404528688960 i n^5 - 1149495367680 i^2 n^5 - 1798761000960 i^3 n^5 
- 1666196421120 i^4 n^5 - 924952089600 i^5 n^5 \nonumber\\
& - 286634557440 i^6 n^5 - 38320128000 i^7 n^5 
- 18618732528 n^6 - 87234630912 i n^6 - 246608991552 i^2 n^6 \nonumber\\
& - 349341189120 i^3 n^6 
- 305355536640 i^4 n^6 - 169087564800 i^5 n^6 - 54286848000 i^6 n^6 - 7664025600 i^7 n^6 
\nonumber\\
& + 13184713168 n^7 + 22506296736 i n^7 + 64354656960 i^2 n^7 + 94757731200 i^3 n^7 
+ 69195962880 i^4 n^7 \nonumber\\
& + 25482885120 i^5 n^7 + 3832012800 i^6 n^7 + 3778269000 n^8 + 11972980536 i n^8 
  + 35129924280 i^2 n^8 \nonumber\\
& + 38594413440 i^3 n^8 + 22345424640 i^4 n^8 + 7025356800 i^5 n^8 + 958003200 i^6 n^8 
- 2732570786 n^9 \nonumber\\
& + 2290135980 i n^9 + 2863365120 i^2 n^9 + 112337280 i^3 n^9 
- 587220480 i^4 n^9 - 159667200 i^5 n^9 \nonumber\\
& - 185879199 n^{10} + 114057636 i n^{10} - 1804196064 i^2 n^{10} - 1755008640 i^3 n^{10} 
- 662618880 i^4 n^{10} \nonumber\\
& - 95800320 i^5 n^{10} 
+ 555823886 n^{11} - 353833392 i n^{11} - 485020800 i^2 n^{11} - 157449600 i^3 n^{11} \nonumber\\
& - 15966720 i^4 n^{11} 
- 37142193 n^{12} - 110297088 i n^{12} + 38272080 i^2 n^{12} + 41247360 i^3 n^{12} 
+ 7983360 i^4 n^{12} \nonumber\\
& - 65980420 n^{13} + 23650440 i n^{13} + 23063040 i^2 n^{13} + 3991680 i^3 n^{13} 
+ 9237162 n^{14} + 7902840 i n^{14} \nonumber\\
& - 807840 i^2 n^{14} - 570240 i^3 n^{14} + 3778940 n^{15} - 1295184 i n^{15} - 475200 i^2 n^{15} 
- 860970 n^{16} - 178200 i n^{16} \nonumber\\
& + 35640 i^2 n^{16} - 48378 n^{17} + 41580 i n^{17} + 29205 n^{18} - 1980 i n^{18} 
- 2970 n^{19} + 99 n^{20}.
\end{align}

%
%
\section{Further Generalizations}

%
%
\subsection{Generalizations to Inverse Factorial Series}
\label{secB1}

Another obvious generalization of our algorithm
would concern the replacement of the 
abscissas by the terms of an inverse factorial
series as opposed to a power series,
\begin{equation}
\label{inverse_factorial}
f(x_i) = \sum_{j=0}^{\infty} c'_j \, x'_{ij}
\approx \sum_{j=0}^{n} c'_j \, x'_{ij}
= c'_0 + \frac{c'_1}{(i+1)} + \frac{c'_2}{(i+1) \, (i+2)} +
\ldots + \frac{c'_n}{(i+1) (i+2) \dots (i+n)} \,,
\end{equation}
where the $x'_{ij}$ are given as follows,
\begin{equation}
\label{absF}
x'_{ij} = \lim_{\epsilon \to 0} \frac{i+\epsilon}{(i+\epsilon)_{j+1}} =
\left\{ \begin{array}{cc}
\dfrac{1}{\Gamma(j+1)} & \qquad \mbox{($i = 0$)} \\[2ex]
\dfrac{i}{(i)_{j+1}} = \dfrac{1}{(i+1) (i+2) \dots (i+j)}  & \qquad \mbox{($i \neq 0$)} 
\end{array} \right. \,.
\end{equation}
We have used the Pochhammer symbol
\begin{equation}
(i)_{j+1} = \frac{\Gamma(i+j+1)}{\Gamma(i)} = i (i+1) \dots (i+j) \,.
\end{equation}
We thus replace the power-like remainder term
estimate~\eqref{remainder_poly} by a remainder estimate
based on an inverse factorial 
series,
\begin{equation}
r'_n = \frac{c'_1}{n+1} 
+ \frac{c'_2}{(n+1) (n+2)} 
+ \frac{c'_3}{(n+1) (n+2) (n+3)} + \dots \,.
\end{equation}
This, in turn, corresponds to Eq.~(8.1-6) of Ref.~\cite{We1989}.
Given $n+1$ partial sums $s_i$ of the input series
$\{ s_i \}_{n=0}^\infty$, the generalized algorithm 
aims to solve the system of equations
\begin{equation}
s_i = \sum_{j=0}^n c'_j(n) \, x'_{ij}
\qquad i, j \in (0, \ldots, n) \,,
\end{equation}
In full analogy to Eq.~\eqref{conv_poly},
one then assumes that the
value of the infinite series is recovered as
$n$ is increased,
\begin{equation}
\label{conv_factorial}
s = s_\infty = \lim_{n \to \infty} c'_0(n) \,.
\end{equation}
We have performed numerical experiments
with the modified Neville scheme based on inverse
factorial series, notably, to the series 
investigated in Secs.~\ref{sec4A} and~\eqref{sec4B}.
This is based on the observation that, if the remainder
term can be expanded in a series with inverse integer powers
of $j$, it can also be expanded into an inverse 
factorial series. In general, the performance of 
the modified algorithm based on inverse factorial series
has been observed to 
be inferior to the inverse-integer-power 
approach outlined in Sec.~\ref{sec2}.
For example, for the input series discussed in Sec.~\ref{sec4A},
one obtains just $7$-decimal convergence for the inverse-factorial
generalized transformation outlined here,  for $n=50$,
while the transformation outlined in Sec.~\ref{sec2} 
leads to $27$-decimal convergence in the same order 
of the transformation. However, one should keep in 
mind that nonlinear sequence transformations
based on remainder estimates with inverse factorials
have shown very favorable numerical results 
for the resummation of divergent 
series~\cite{Si1979,Si1980,Si1986,We1989,We1996c,Si2003}.
The relative performance of the generalized enhanced Neville scheme based on
inverse factorial series as compared to inverse power series
is expected to significantly depend on the problem
under study, and the method outlined 
in Eq.~\eqref{inverse_factorial}---\eqref{conv_factorial}
could very well become useful in other contexts.

%
%
\subsection{Generalizations to Quadratic and Quartic Remainder Terms}
\label{secB2}

In this section we present a further generalization of the enhanced Neville
algorithm adapted to series involving second and fourth powers:
\begin{equation}
f(x_i) = \sum_{j=0}^{\infty} c^{(\alpha)}_j \, 
x^{(\alpha)}_{ij}
\approx 
\sum_{j=0}^{n} c^{(\alpha)}_j \, x^{(\alpha)}_{ij}
= c^{(\alpha)}_0 + 
\frac{c^{(\alpha)}_1}{(i+1)^\alpha} + 
\frac{c^{(\alpha)}_2}{(i+1)^{2\alpha}} +
\ldots + \frac{c^{(\alpha)}_n}{(i+1)^{n\alpha}}, \;\;\; \alpha=2,4
\label{eq:polyAlpha}
\end{equation}
where the terms $x^\alpha_{ij}$ are defined as,
\begin{equation}
\label{absAlpha}
x^{(\alpha)}_{ij} = \frac{1}{(i+1)^{j\alpha}}.
\end{equation}
Our modified remainder term
is thus based on second and fourth powers,
\begin{equation}
r^{(\alpha)}_n 
= \frac{c^{(\alpha)}_1}{(i+1)^{\alpha}} 
+ \frac{c^{(\alpha)}_2}{(i+1)^{2\alpha}} 
+ \frac{c^{(\alpha)}_3}{(i+1)^{3\alpha}} + \dots \,,
\qquad \alpha=2,4.
\end{equation}

In these two cases, one may still derive an analytic 
expression for the coefficient $c^{(\alpha)}_0(n)$, 
which approaches the value of the infinite series as $n$ approaches infinity:
\begin{equation}
\label{conv_alpha}
s = s_\infty = \lim_{n \to \infty} c^{(\alpha)}_0(n) \,.
\end{equation}
The coefficients can be explicitly calculated as follows:
\begin{equation}
\begin{split}
c^{(2)}_0(n) &= \sum_{i=0}^n 
\frac{(-1)^{n-i} 2 (i+1)^{2n+2}}{\Gamma(i+n+3) \Gamma(n-i+1)}s_i,\\
c^{(4)}_0(n) &= \sum_{i=0}^n \frac{\pi \, 
\text{csch}((i+1)\pi) \, (-1)^{n-i} \, 4 \, (i+1)^{1+4(n+1)}}%
{\Gamma(n-i+1) \,
\Gamma(i+n+3) \, \Gamma(-(i+1) \, \textrm{j} + n + 2) \,
\Gamma((i+1) \, \textrm{j} + n + 2)} \, s_i \,,
\end{split}
\end{equation}
where $\textrm{j}$ is the imaginary unit. (The notation
$\textrm{j}$ is used in order to clearly distinguish 
the imaginary unit from the summation index $i$.) We 
have performed numerical experiments for 
the input series discussed in this 
article, and for other input series, found that the algorithms with
$\alpha = 2$ and $\alpha=4$ typically 
perform well only in rare cases,
in which the remainder term involves only 
even integer powers of $1/(n+1)$. 
Such cases are rare. For example, in the 
conceptually simple case of
\begin{equation}
s_n = \sum_{k=0}^n \frac{1}{k^3} \,,
\qquad
s_\infty = \zeta(3) \,,
\qquad
r_n = -\sum_{k=n+1}^\infty \frac{1}{k^3} = \psi^{(2)}(n+1) 
= -\frac{1}{2 n^2} + \frac{1}{2 n^3} - \frac{1}{4 n^4} + \dots \,,
\end{equation}
the remainder term still possess an expansion,
for large $n$, containing odd powers of $n$,
even if the input series contains only odd powers
of $k$. Its remainder term, under the assumption that 
summation could be replaced by integration,
would be assumed to contain only even powers 
of $n$, but this is not the case. Here, $\psi^{(m)}(z) = (\dd^{m+1}/\dd z^{m+1})
\ln \Gamma(z)$ is the $(m+1)$th logarithmic derivative 
of the Gamma function. The presence of odd powers in 
the remainder term persists if the 
remainder term is written as a function of $n+1$ as opposed to $n$.
Still, for completeness, we believe it is useful to indicate 
the generalization of the algorithm to series 
whose remainder terms contain only even inverse powers of $n+1$.

%
%
\section{Bethe Logarithms for Excited States}
\label{appb}

We supplement the results presented
in Sec.~\ref{sec4} with a somewhat more detailed
discussion.
The Lerch $\Phi$ transcendent~\cite{AbSt1972} is defined as
\begin{equation}
\Phi(z,s,a) = \sum_{n=0}^\infty \frac{z^n}{(n+a)^s} \,.
\end{equation}
A formula valid in a larger domain of the 
complex plane is 
\begin{equation}
\Phi(z,s,a) =  
\frac{1}{1-z} \; \sum_{n=0}^\infty \left( \frac{-z}{1-z} \right)^n 
\sum_{k=0}^n \frac{(-1)^k}{(a+k)^{s}} \, {n \choose k} \,.
\end{equation}
The Lerch transcendent can be thought of as a generalization
of the zeta function, because $\Phi(1,s,1) = \zeta(s)$.
The Lerch transcendent, for the case $s=1$, 
can be written as a hypergeometric function, 
$\Phi(z,1,a) = {}_2 F_1(1,a,1+a,z)/a$.
According to Chap.~4 of Ref.~\cite{JeAd2022book}, hydrogenic 
Bethe logarithms of $1S$, $2S$ and $2P$ 
states can be written as rather compact sums over Lerch $\Phi$ 
functions, and evaluated to essentially arbitrary accuracy
using the convergence acceleration algorithms outlined here.
We have for the ground state,
\begin{align}
\ln k_0(1S) = & \;
10 \, \ln(2) - 2 \zeta(2) - 1 + \sum_{k=2}^\infty 
\frac{16 \, k}{(k-1)^2 \, (k+1)^2} \,
\Phi\left(\frac{1+k}{1-k}, 1, 2 \, k\right)
\nonumber\\[2ex]
=& \;
2.%
98412\,85557\,65497\,61075\,97770\,90013\,79796\,99751\,80566\,17002\,%
\nonumber\\
& \,\;\;\;\,00048\,15926\,13924\,06576\,62306\,75532\,86860\,62013\,30404\,72249 \,.
\end{align}
This result adds 50~figures to the result communicated
in Eq.~\eqref{lnk01S}.
The logarithmic sum, for the $2S$ state, is given by
\begin{align}
\ln k_0(2S) = & \;
-\frac{545}{36} + \frac{16}{3} \,\ln(2) - 14 \, \zeta(2) + 24 \, \zeta(3)
+ \sum_{k=3}^\infty 
\frac{1024 \, k \, (k-1)\,(k+1)}{(k-2)^3 \, (k+1)^3} \,
\Phi\left(\frac{2+k}{2-k}, 1, 2 \, k\right)
\nonumber\\[2ex]
=& \;
2.%
81176\,98931\,20563\,51521\,97427\,85941\,63611\,28935\,51470\,29732%
\nonumber\\
& \,\;\;\;\,41909\,18696\,96453\,24020\,20118\,89106\,87017\,48612\,02831\,24031 \,,
\end{align}
whereas for the $2P$ state, it reads
\begin{align}
\ln k_0(2P) = & \;
-\frac{3437}{2916} + \frac{3280}{2187} \,\ln(2) - \frac{14}{3} \, \zeta(2) 
+ \frac{136}{9} \, \zeta(3) - \frac{64}{3} \, \zeta(4)
\nonumber\\[2ex]
& \; + \sum_{k=3}^\infty 
\frac{256\, k^3 \, (11\,k^2-12)}{3 \, (k-2)^4 \, (k+2)^4} \,
\Phi\left(\frac{2+k}{2-k}, 1, 2 \, k\right)
\nonumber\\[2ex]
=& \;
-0.%
03001\,67086\,30212\,90244\,36757\,10951\,14406\,39409\,33044\,23103
\nonumber\\
& \,\;\;\;\;\;\;\;\,%
04668\,98525\,32719\,44796\,89622\,57183\,26244\,10312\,70799\,73828 \,.
\end{align}
While these formulas do not provide direct, analytic results 
for Bethe logarithms, they still illustrate that one 
can think of the Bethe logarithm as a sum of logarithms, 
zeta functions, and generalized zeta functions.


\begin{thebibliography}{10}

\bibitem{JeMoSo1999}
U.~D. Jentschura, P.~J. Mohr, and G. Soff, {\em \relax{Calculation of the
  Electron Self-Energy for Low Nuclear Charge}},  Phys. Rev. Lett. {\bf 82},
  53--56  (1999).

\bibitem{JeMoSoWe1999}
U.~D. Jentschura, P.~J. Mohr, G. Soff, and E.~J. Weniger, {\em
  \relax{Convergence acceleration via combined nonlinear-condensation
  transformations}},  Comput. Phys. Commun. {\bf 116},  28--54  (1999).

\bibitem{JeAd2022book}
U.~D. Jentschura and G.~S. Adkins, {\em \relax{Quantum Electrodynamics: Atoms,
  Lasers and Gravity}} (World Scientific, Singapore, 2022).

\bibitem{GiEtAl2020}
L.~T. Giorgini, U.~D. Jentschura, E.~M. Malatesta, G. Parisi, T. Rizzo, and J.
  Zinn-Justin, {\em \relax{Two--Loop Corrections to the Large--Order Behavior
  of Correlation Functions in the One--Dimensional $N$--Vector Model}},  Phys.
  Rev. D {\bf 101},  125001  (2020).

\bibitem{GiEtAl2022}
L.~T. Giorgini, U.~D. Jentschura, E.~M. Malatesta, G. Parisi, T. Rizzo, and J.
  Zinn-Justin, {\em \relax{Correlation Functions of the Anharmonic Oscillator:
  Numerical Verification of Two–Loop Corrections to the Large–Order
  Behavior}},  Phys. Rev. D {\bf 105},  105012  (2022).

\bibitem{GiEtAl2024i}
L.~T. Giorgini, U.~D. Jentschura, E.~M. Malatesta, G. Parisi, T. Rizzo, and J.
  Zinn-Justin, {\em Instantons in $\phi^4$ Theories: Transseries, Virial
  Theorems and Numerical Aspects (e-print arXiv:2405.18191 [hep-th])},
  to appear in Physical Review D (2024).

\bibitem{We1989}
E.~J. Weniger, {\em Nonlinear sequence transformations for the acceleration of
  convergence and the summation of divergent series},  Comput. Phys. Rep. {\bf
  10},  189--371  (1989).

\bibitem{BrRZ1991}
C. Brezinski and M. Redivo-Zaglia, {\em \relax{Extrapolation Methods}}
  (North-Holland, Amsterdam, 1991).

\bibitem{CaEtAl2007}
E. Caliceti, M. Meyer-Hermann, P. Ribeca, A. Surzhykov, and U.~D. Jentschura,
  {\em \relax{From Useful Algorithms for Slowly Convergent Series to Physical
  Predictions Based on Divergent Perturbative Expansions}},  Phys. Rep. {\bf
  446},  1--96  (2007).

\bibitem{Ai1926}
A.~C. Aitken, {\em \relax{On {B}ernoulli's numerical solution of algebraic
  equations}},  Proc. Roy. Soc. Edinburgh {\bf 46},  289--305  (1926).

\bibitem{Sh1955}
D. Shanks, {\em Non-linear transformations of divergent and slowly convergent
  sequences},  J. Math. and Phys. (Cambridge, Mass.) {\bf 34},  1--42  (1955).

\bibitem{Wy1956eps}
P. Wynn, {\em \relax{On a device for computing the $e_m (S_n)$
  transformation}},  Math. Tables Aids Comput. {\bf 10},  91--96  (1956).

\bibitem{Je2000prd}
U.~D. Jentschura, {\em \relax{Resummation of nonalternating divergent
  perturbative expansions}},  Phys. Rev. D {\bf 62},  076001  (2000).

\bibitem{HeEu1936}
W. Heisenberg and H. Euler, {\em \relax{Folgerungen aus der Diracschen Theorie
  des Positrons}},  Z. Phys. {\bf 98},  714--732  (1936).

\bibitem{Sc1951}
J. Schwinger, {\em \relax{On Gauge Invariance and Vacuum Polarization}},  Phys.
  Rev. {\bf 82},  664--679  (1951).

\bibitem{DiRe1985}
W. Dittrich and M. Reuter, {\em Effective Lagrangians in Quantum
  Electrodynamics -- Lecture Notes in Physics Vol. 220} (Springer, Berlin,
  Heidelberg, New York, 1985).

\bibitem{Si1979}
A. Sidi, {\em \relax{Convergence properties of some nonlinear sequence
  transformations}},  Math. Comput. {\bf 33},  315--326  (1979).

\bibitem{Si1980}
A. Sidi, {\em \relax{Analysis of the convergence of the $T$-transformation for
  power series}},  Math. Comput. {\bf 35},  833--850  (1980).

\bibitem{Si1986}
A. Sidi, {\em \relax{Borel summability and converging factors for some
  everywhere divergent series}},  SIAM J. Math. Anal. {\bf 17},  1222--1232
  (1986).

\bibitem{We1996c}
E.~J. Weniger, {\em Computation of the Whittaker function of the second kind by
  summing its divergent asymptotic series with the help of nonlinear sequence
  transformations},  Comput. Phys. {\bf 10},  496--503  (1996).

\bibitem{Si2003}
A. Sidi, {\em Practical Extrapolation Methods} (Cambridge University Press,
  Cambridge, 2003).

\bibitem{NP1961}
C.~W. Clenshaw, E.~T. Goodwin, D.~W. Martin, G.~F. Miller, F.~W.~J. Olver, and
  J.~H. Wilkinson, {\em \relax{(National Physical Laboratory), Modern Computing
  Methods, Notes on Applied Science}}, 2 ed. (H. M. Stationary Office, London,
  1961), Vol.~16.

\bibitem{AkSaJeBeSoMo2003}
S.~V. Aksenov, M.~A. Savageau, U.~D. Jentschura, J. Becher, G. Soff, and P.~J.
  Mohr, {\em \relax{Application of the combined nonlinear-condensation
  transformation to problems in statistical analysis and theoretical physics}},
   Comput. Phys. Commun. {\bf 150},  1--20  (2003).

\bibitem{Br1980E}
C. Brezinski, {\em \relax{A general extrapolation algorithm}},  Numer. Math.
  {\bf 35},  175--180  (1980).

\bibitem{Ha1979E}
T. Havie, {\em \relax{Generalized {N}eville type extrapolation schemes}},
  B.I.T. {\bf 19},  204--213  (1979).

\bibitem{Ne1934}
E.~H. Neville, {\em \relax{Iterative Interpolation}},  J. Indian Math. Soc.
  {\bf 20},  87--120  (1934).

\bibitem{GaGu1974}
D.~S. Gaunt and A.~J. Guttmann,  in {\em \relax{Phase Transitions and Critical
  Phenomena 3}}, edited by C. Domb and M.~S. Green (Academic Press, London,
  1974), p.\ 181.

\bibitem{vT1994}
A.~H. van Tuyl, {\em \relax{Acceleration of Convergence of a Family of
  Logarithmically Convergent Series}},  Math. Comput. {\bf 63},  229--246
  (1994).

\bibitem{Be1947}
H.~A. Bethe, {\em \relax{The Electromagnetic Shift of Energy Levels}},  Phys.
  Rev. {\bf 72},  339--341  (1947).

\bibitem{FeBa1992}
H.~R.~P. Ferguson and D.~H. Bailey, {\em A Polynomial Time, Numerically Stable
  Integer Relation Algorithm}, RNR Techn. Rept. RNR--91--032 (1992).

\bibitem{BaPl1997}
D.~H. Bailey and S. Plouffe,  in {\em \relax{Organic Mathematics: Proceedings
  of the Workshop Held in Burnaby, BC}}, edited by J. Borwein, P. Borwein, L.
  J\"{o}rgenson, and R. Corless (American Mathematical Society, Philadelphia,
  PA, 1997), pp.\ 73--88.

\bibitem{FeBaAr1999}
H.~R.~P. Ferguson, D.~H. Bailey, and S. Arno, {\em \relax{Analysis of PSLQ, an
  integer relation finding algorithm}},  Math. Comput. {\bf 68},  351--369
  (1999).

\bibitem{BaBr2001}
D.~H. Bailey and D. Broadhurst, {\em \relax{Integer Relation Detection}},
  Math. Comput. {\bf 70},  1719--1736  (2001).

\bibitem{PrFlTeVe1993}
W.~H. Press, B.~P. Flannery, S.~A. Teukolsky, and W.~T. Vetterling, {\em
  \relax{Numerical Recipes in {C}: The Art of Scientific Computing}}, 2 ed.
  (Cambridge University Press, Cambridge, UK, 1993).

\bibitem{AbSt1972}
M. Abramowitz and I.~A. Stegun, {\em \relax{Handbook of Mathematical
  Functions}}, 10 ed. (National Bureau of Standards, Washington, D. C., 1972).

\end{thebibliography}
\end{document}